\documentclass[letterpaper, 10 pt, conference]{ieeeconf}

\pdfoutput=1

\usepackage[vince]{irom}

\usepackage{multicol}
\usepackage[caption=false]{subfig}
\usepackage{csquotes}
\usepackage[american]{babel}
\usepackage{dblfloatfix}
\usepackage{ifthen}
\usepackage{placeins}
\usepackage{setspace}
\usepackage{titlesec}

\usepackage[backend=bibtex,
    citestyle=numeric-comp,
    bibstyle=ieee,
    maxnames=2,
    minnames=2,
    doi=false,
    isbn=false,
    url=false,
    eprint=false]{biblatex}

\usepackage[protrusion=true,
    activate={true,nocompatibility},
    final,
    tracking=true,
    kerning=true,
    spacing=true,
    factor=1100]{microtype}
    
\usepackage{siunitx}

\newcommand{\OPToff}{$\mathrm{OFF}$}
\newcommand{\OPTon}{$\mathrm{ON}$}

\newcommand{\off}[1]{{#1}^{\mathrm{off}}}

\newcommand{\on}[1]{{#1}^{\mathrm{on}}}
\newcommand{\ol}[1]{{#1}^{\perp}}
\newcommand{\opt}[1]{{#1}^{\mathrm{opt}}}

\newcommand{\func}[1]{\mathds{#1}}
\newcommand{\mech}[1]{\mathcal{#1}}
\newcommand{\set}[1]{\mathbf{#1}}

\newcommand{\statespace}{\set{X}}
\newcommand{\inputspace}{\set{U}}

\makeatletter
\newcounter{IEEE@bibentries}
\renewcommand\IEEEtriggeratref[1]{%
  \renewbibmacro{finentry}{%
    \stepcounter{IEEE@bibentries}%
    \ifthenelse{\equal{\value{IEEE@bibentries}}{#1}}
    {\finentry\@IEEEtriggercmd}
    {\finentry}%
  }%
}
\makeatother

\declaretheorem[name=Proposition, style=definition]{proposition}

\declaretheorem[name=Remark, style=definition]{remark}
\declaretheorem[numbered=no, name=Definition, style=definition]{definition}

\newtoggle{includeappendix}
\togglefalse{includeappendix}

\SetTracking{encoding={*}, shape=sc}{10}

\titlespacing{\section}{0pt}{2ex}{1ex}
\titlespacing{\subsection}{0pt}{1ex}{1ex}
\def\baselinestretch{0.99}

\IEEEoverridecommandlockouts
\hypersetup{pdfauthor={Vincent Pacelli}, pdftitle={Robust Control Under Uncertainty via Bounded Rationality and Differential Privacy}}

\begin{document}

\title{\LARGE \bf Robust Control Under Uncertainty via Bounded Rationality \\ and Differential Privacy}

\author{Vincent Pacelli and Anirudha Majumdar
    \thanks{The authors are affiliated with the Mechanical and Aerospace Engineering department of Princeton University, Princeton, NJ, 08544, USA {\tt\footnotesize \{vpacelli, ani.majumdar\}@princeton.edu}.}
    \thanks{This work is partially supported by the Office of Naval Research [N00014-21-1-2803] and the NSF CAREER award [2044149].}
}

\maketitle

\begin{abstract}
The rapid development of affordable and compact high-fidelity sensors (e.g., cameras and LIDAR) allows robots to construct detailed estimates of their states and environments. However, the availability of such rich sensor information introduces two technical challenges: (i) the lack of analytic sensing models, which makes it difficult to design controllers that are robust to sensor failures, and (ii) the computational expense of processing the high-dimensional sensor information in real time. This paper addresses these challenges using the theory of \emph{differential privacy}, which allows us to (i) design controllers with bounded sensitivity to errors in state estimates, and (ii) bound the amount of state information used for control (i.e., to impose \emph{bounded rationality}). The resulting framework approximates the \emph{separation principle} and allows us to derive an upper-bound on the cost incurred with a faulty state estimator in terms of three quantities: the cost incurred using a perfect state estimator, the magnitude of state estimation errors, and the level of differential privacy. We demonstrate the efficacy of our framework numerically on different robotics problems, including nonlinear system stabilization and motion planning. 

\end{abstract}

\IEEEpeerreviewmaketitle

\section{Introduction}
\label{sec:intro}

Despite the increasing availability of high-resolution sensors for robotic systems, \emph{partial-observability} remains a challenge for controlling such systems. Sensing modalities such as vision and LIDAR are ultimately noisy and only provide partial information about the robot's state and environment. In general, solving optimal control problems with partial observability is computationally intractable. One of the most common approaches to tackling this challenge is to assume the \emph{separation principle} \cite{Anderson07}, i.e., to \emph{independently} design (i) a state estimator, and (ii) a controller, e.g., one based on model-predictive control (MPC), that is optimal assuming perfect state estimation. The modularity afforded by such an approach coupled with the relative tractability of tackling the estimation and control problems independently make this framework appealing. However, the separation principle does not hold in general for robotic systems due to nonlinear dynamics / measurement models. A controller that assumes perfect state estimation can thus be highly sensitive to small errors in the state estimate, leading to significant brittleness of the overall control system. This is particularly challenging in the increasingly common case where a (deep) learning model is used as part of the robot's state estimation pipeline (due to potential over-fitting). As a result, robots often behave erratically when faced with unforeseen measurement errors despite possessing high-resolution sensors. 

In contrast, the cognitive science literature on \emph{bounded rationality} demonstrates that humans display impressive levels of robustness and generalization in dexterous tasks such as locomotion and ball-catching without relying on highly-accurate state estimates \cite{Gigerenzer09, Gigerenzer11}. For example, \emph{gaze heuristics} are control laws used by humans to catch freely-falling objects by adjusting their running speed based on the motion of the object in their visual field \cite{Hofer18, Belousov16}. Unlike MPC strategies, gaze heuristics do not require an accurate estimate of the system state or other quantities such as the wind speed or object mass. In addition to the robustness afforded by such heuristics, they are often extremely computationally efficient. 
Interestingly, the dual benefits of robustness and computational efficiency are related to each other in the above heuristics --- boundedness of online computation (i.e., bounded rationality) ensures that control actions are only loosely coupled to sensor measurements. Such a loose coupling can prevent uncertainty in measurements from significantly impacting controller performance. 

\begin{figure}[!t]
    \centering
    \begin{minipage}[t]{\linewidth}
        \centering
        {\fontsize{11}{18} \textrm{Motion Planning with Sensing Uncertainty}}\par
        \vspace{2pt}
        \includegraphics[width=0.49\linewidth]{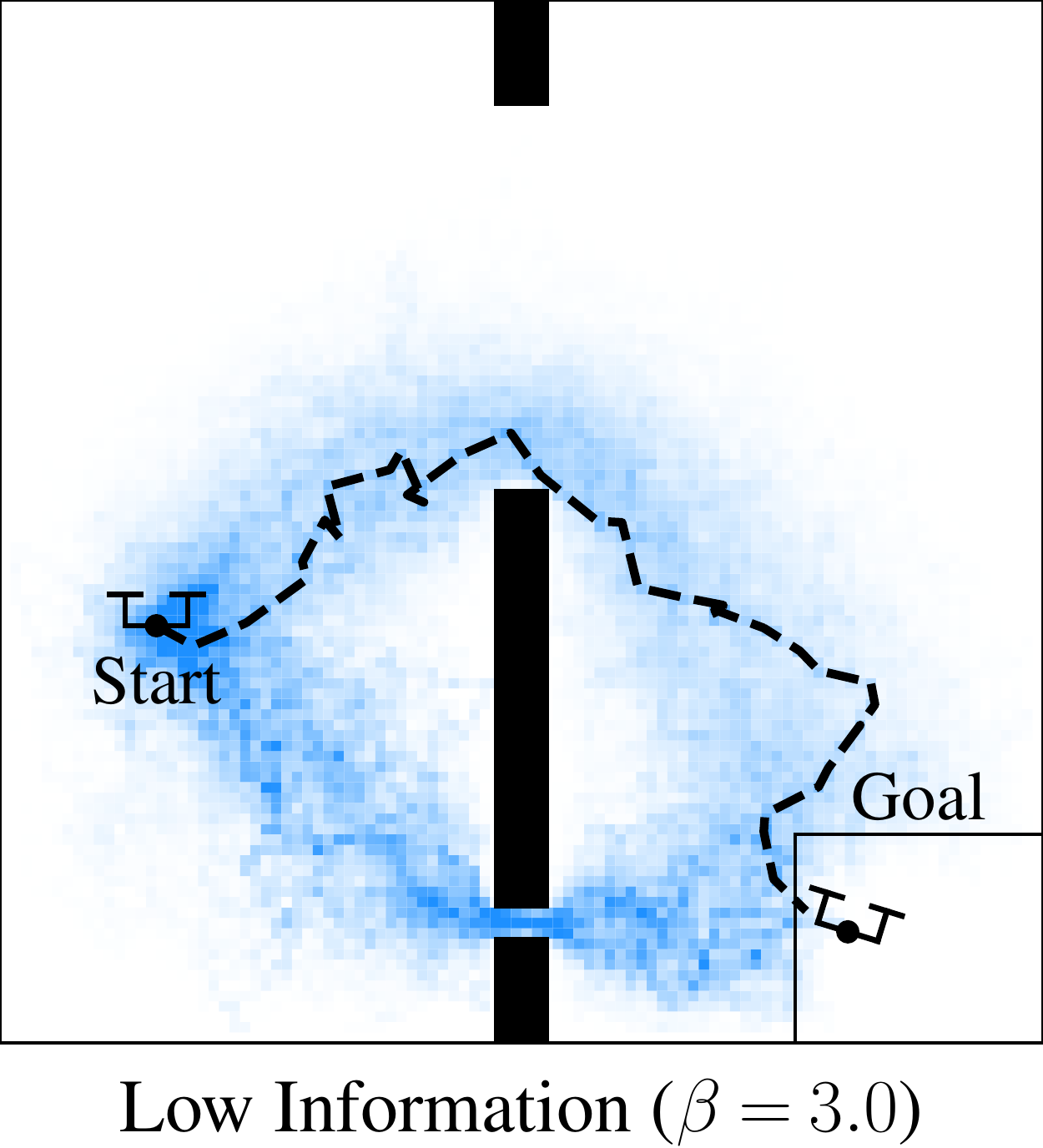}
        \hfill
        \includegraphics[width=0.49\linewidth]{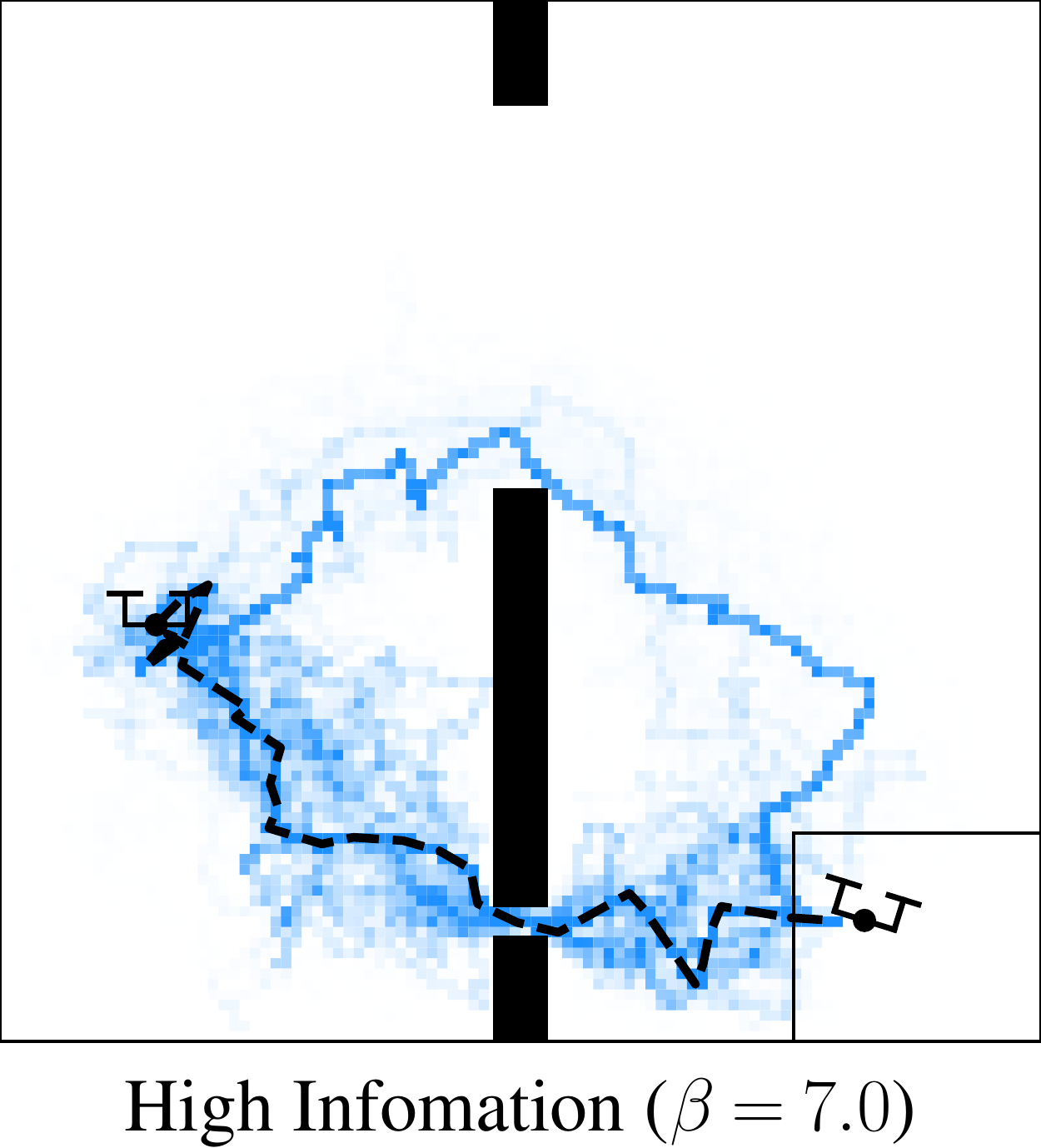}
    \end{minipage}
    
    \vspace{8pt}
    \includegraphics[width=\linewidth]{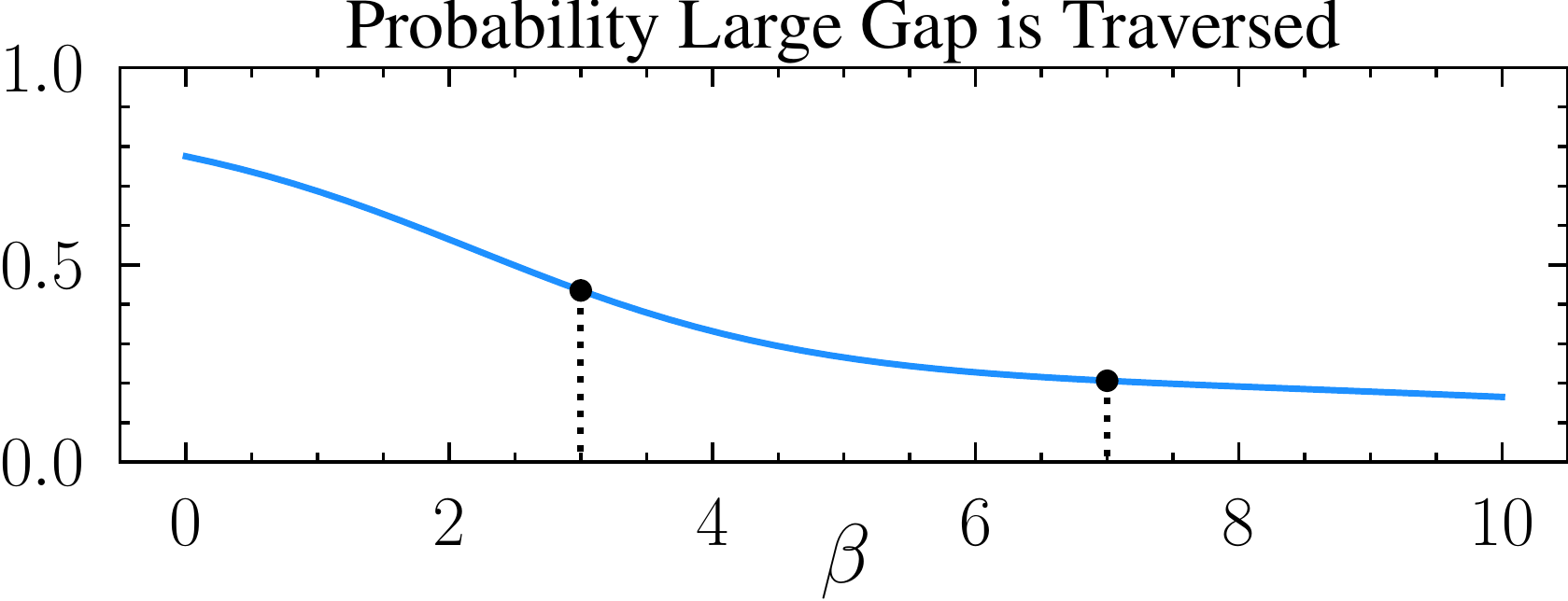}
    
    \vspace{-8pt}
    \caption{The information-constrained motion-planning problem  (\cref{sec:examples}). The robot must navigate from the start to the goal using a limited amount of state information (correlated with $\beta$). The distribution of trajectories is visualized for two values of $\beta$. At low information, the robot prefers the longer path through the wider gap, which is more robust to positional uncertainty than the direct route preferred at higher information usage.}
    \label{fig:anchor}
    \vspace{-1mm}
\end{figure}

\textbf{Statement of Contributions.} The key conceptual contribution of this work is to use the observation that bounded rationality is linked to robust decision-making as a guide and formalize it using the theory of \emph{differential privacy} (DP) \cite{Dwork06, Dwork11, Dwork14}. Originally developed as a framework for ensuring the privacy of individuals’ data (e.g., in a database while publicly releasing statistics), DP allows us to formalize the idea that control inputs should not depend too tightly on state estimates. In particular, our approach preserves the modular structure of separation principle-based approaches by crucially ensuring that the mapping from state estimates to control inputs is differentially private. Moreover, it exploits the interpretation of the differentially private \emph{exponential mechanism} as synthesizing \emph{information-constrained} (bounded-rational) controllers.

The primary theoretical contribution of this paper is a novel bound on the expected performance of such controllers in the presence of estimation error. The bound depends only on the cost incurred using a perfect state estimator, the magnitude of estimation errors, and the level of differential privacy.
Moreover, due to the choice of the differentially-private mechanism on which our policy is based, our bound naturally applies to many popular inference-based control algorithms, such as optimal control methods based on importance sampling, Stein-variational gradient descent, and $Q$-learning. We are thus able to demonstrate the efficacy of the bounded-rationality approach we propose on a number of robotics problems of interest, such as stabilizing a nonlinear planar quadrotor and robust motion planning. To our knowledge, the approach presented in this paper is the first to utilize the framework of differential privacy to achieve robust and bounded-rational control of robotic systems.

\section{Related Work}
\label{sec:related}

This section briefly reviews three areas of the literature that are typically considered independently from one another. However, the unifying theme is that they exploit various properties of the Gibbs measure.\footnote{The Gibbs measure is also known as the exponential family, softargmax, or Boltzmann distribution in different contexts and fields.} This work will exploit developments in each cited field.

\subsection{Differential Privacy and Applications in Control Theory}

Differential privacy (DP)  \cite{Dwork06, Dwork11, Dwork14, Chatzikokolakis13} is an algorithmic framework that arose to meet the conflicting needs of statisticians (who ultimately want to publish analyses of sensitive data sets) and participants (who want to keep their data private). 
DP offers an elegant solution based on the intuition that if the results published from a study are insensitive to the substitution of data from any given individual, then this offers privacy. The literature on DP is vast and many differentially-private \emph{mechanisms} (i.e., stochastic algorithms) have been proposed. The most relevant for our purpose is the \emph{exponential mechanism} \cite{Dwork14} defined by the Gibbs measure; this mechanism provides a method for releasing the results of an optimal decision-making procedure based on data (e.g., Bayesian inference \cite{Dimitrakakis17}). Additionally, DP enjoys several appealing theoretical properties, namely the preservation of privacy under composition of mechanisms, and the ability to rigorously analyze the trade-off between privacy and accuracy \cite{Dwork14}.

In the context of control theory, application of the DP formalism is largely limited to networked control and distributed systems \cite{Le13, Cortes16, Nozari16, Han16, Koufogiannis17, Han18}. These include privately solving distributed optimization problems, modifying control inputs to keep the system state private from an external observer, and aggregating measurements from privacy-concerned agents for filtering performed by a central entity. In contrast, we \emph{do not} use DP for the sake of privacy; instead, DP formalizes the intuition described in \cref{sec:intro}: limiting the sensitivity of control inputs to state estimates affords robustness to estimation error. To our knowledge, our work is the first to utilize DP for this purpose.

\subsection{Bounded Rationality and Robust Decision-Making}

Bounded rationality is a model of decision-making first introduced to address misalignment between economic theory based on rational agents and the reality of sub-optimal human decision-making \cite{Simon55}. In short, while agents often have an objective they are trying to optimize, their ability to optimize this objective is bounded by informational and computational constraints. The cognitive science literature has identified several \emph{heuristics} used by humans to deal with such constraints \cite{Warren86, Lee76, Gigerenzer09, Gigerenzer11}. Empirical work in cognitive science and robotics suggests that in addition to efficiency, such heuristics can provide generalization to new environments \cite{Gigerenzer09, Gigerenzer11, Hofer18, Pacelli19, Pacelli20}. 

While bounded rationality is formalized in a number of ways for artificial agents, the most relevant framework adds an \emph{information-theoretic constraint} to an otherwise rational agent. More precisely, a bounded-rational agent is modeled as solving an \emph{entropy-maximization} or \emph{rate-distortion} problem \cite{Tishby11, Fox16, Fox12, Ortega14, Ortega15, Shafieepoorfard16}. Despite empirical evidence that such a bounded-rational agent can be robust to uncertainty or noise in sensor measurements and dynamics \cite{Pacelli19, Igl19, Pacelli20, Lu20, Pedram21}, formal analyses justifying such robustness benefits are limited. One approach is to use a variational representation of the (relative) entropy (e.g., the Donsker-Varadhan representation \cite{Donsker83}) to derive a zero-sum game that the bounded-rational agent is implicitly playing against an adversary that chooses a cost function \cite{Ortega14, Ortega15, Eysenbach21}. This approach is closely related to the \emph{maximum-entropy principle} from Bayesian statistics \cite{Grunwald04} and leads to generalization results for some estimation problems. However, the payoff of the game is difficult to interpret (especially for the adversary) in the context of most control problems, which limits the usefulness of the analysis. Alternatively, an upper bound on the performance degradation of a bounded-rational agent under measurement error exists \cite{Pacelli19} using variational representations of the entropy, but its assumptions make the bound difficult to interpret practically. Instead, this paper demonstrates that the connection between DP and bounded rationality through the Gibbs measure allows for the derivation of a bound on control performance in the presence of estimation error through a large-deviation analysis. The trade-offs presented by this bound are easily understood and, importantly, the bound only depends on access to simulation or lab data without actually deploying the robot in the target environment.


\subsection{Optimal Control as Bayesian Inference}
\label{ssec:optimal control}

Another relevant vein of research at the intersection of information theory, Bayesian inference, and control theory is \emph{linearly-solvable optimal control} (LSOC) \cite{Kappen05, Kappen05b, Todorov08, Todorov09, Williams17, Theodorou12, Braun11, Lambert20, Barcelos21, Horowitz14}. These techniques exploit an identified equivalence between optimal control and Bayesian inference: an approximately optimal control sequence is found by sampling from a Gibbs measure over control inputs conditioned on the current state and the event that the sequence results in an optimal trajectory \cite{Todorov08}. The result is that the nonlinear stochastic optimal control problem is solved as a linear, albeit infinite-dimensional, differential equation using a transformed cost function, and the solution is approximately computed using inference algorithms that include importance sampling \cite{Murphy12}, Stein-variational gradient descent (SVGD) \cite{Liu16}, and $Q$-learning \cite{Haarnoja17}. Despite the growing popularity of these algorithms for optimal control and empirical evidence for the generalization benefits they confer \cite{Haarnoja18}, theoretical knowledge of their robustness properties is limited \cite{Levine18}. This paper aims to fill this gap by providing a new, concrete analysis that these algorithms are provably robust to estimation error.

\section{Notation}
\label{sec:notation}
Random variables are denoted by uppercase letters (e.g., $X$), and realized quantities are denoted by lowercase letters (e.g., $x$). Deterministic functions appear in either case. Finite sequences are represented as $x_{i:j} = (x_k)_{k = i}^{j}$ for $i \leq j$. The indicator function for a set $\set{A}$ is denoted $1_\set{A}(\cdot)$. Functionals are double-struck, namely decorated varieties of the expectation $\func{E}[\cdot]$, the relative entropy $\func{D}[\cdot || \cdot]$, and the tightest Lipschitz constant of a scalar-valued function $\func{L}[\cdot]$. Expectations are taken over the uppercase random variables (and mechanisms), e.g., in $\func{E}[H(x, U)]$, $x$ is fixed and integration is over $U$. Sets are in boldface and $\set{\Delta}(\set{A})$ is the set of distributions with support $\set{A}$. For brevity, it is assumed that all necessary moments of random variables exist, and spaces are measurable with their subsets coming from appropriate $\sigma$-algebras. 

A \emph{mechanism} (denoted by uppercase script) refers to a (randomized) algorithm (formally, a transition kernel) between sets $\set{X}$ and $\set{Y}$. A mechanism $\mech{M}: \set{X} \to \set{\Delta}(\set{Y})$ defines a probability distribution on $\set{Y}$ for each $x \in \set{X}$. Denote by $\mech{M}(x)\{\cdot\}$ the density and measure of this distribution when applied to elements and measurable subsets of $\set{Y}$ respectively. When clear from context, we will overload this notation by treating $\mech{M}(x)$ as a random variable with support $\set{Y}$. Mechanisms may be composed, i.e., if $\mech{M}': \set{Y} \to \set{\Delta}(\set{Z})$ is another mechanism, then $(\mech{M}' \circ \mech{M})(\cdot) \coloneqq \mech{M}'(\mech{M}(\cdot))$.


\section{Robust Single-Step Decision-Making}
\label{sec:robust}
This section details the robustness properties that follow from applying a bounded-rationality approach to a single-step decision-making problem. \cref{sec:multistep} utilizes composition properties of DP to extend the analysis to multi-step optimal control problems.

\subsection{Problem Statement}
\label{ssec:problem statement}

Let the state of the robotic system be $x \in \set{X}$. The goal of the agent (robot) is to process the information contained in this state and select a control input $u \in \set{U}$ that minimizes a cost $H(x, u)$. The state $X \sim \mech{X}$ is random and the agent must find a feedback mechanism $\mech{U}: \set{X} \to \set{\Delta}(\set{X})$ that solves:
\begin{align}
    \tag{\OPToff}
    \min_{\mech{U}} \off{\func{J}}[\mech{U}] \coloneqq \func{E}[H(X, \mech{U}(X))]. \label{eq:optoff}
\end{align}
Since $X$ is made available to the agent, this problem is fully-observable. It is referred to as the \emph{offline problem} since it corresponds to, e.g. a lab or simulation setting where the agent is able to make arbitrarily fine measurements about the world.

In practice, the primary concern is performance of the feedback controller on the \emph{online problem}, where the agent only has access to a noisy state estimate $\hat{X} \sim \hat{\mech{X}}(X)$.\footnote{Typically $\hat{\mech{X}}(x)$ is the composition of a state estimator with a noisy sensor, but these mappings are implementation dependent and irrelevant to the analysis. It is simpler and without loss of generality to work only with their composition.} The control input $U \sim \mech{U}(\hat{X})$ is selected using this estimate, and the goal is to solve:
\begin{align}
\tag{\OPTon}
    \min_{\mech{U}} \on{\func{J}}[\mech{U}] \coloneqq \off{\func{J}}[\mech{U} \circ \hat{\mech{X}}]. \label{eq:opton}
\end{align}
This is a \emph{partially-observable} decision problem and its general solution requires reformulating the problem into an intractably high-dimensional (often infinite-dimensional) optimization problem \cite{Bertsekas17, James96}. Instead, we will demonstrate that a \emph{bounded-rational} controller allows for the value of \cref{eq:opton} to be bounded in terms of \cref{eq:optoff}. That is, \emph{performance on the online problem can be guaranteed using only information available offline} due to a property of the bounded rationality controller known as \emph{differential privacy}. The difference in performance between these two problems is defined to be $\Delta \func{J}[\mech{U}] \coloneqq \on{\func{J}}[\mech{U}] - \off{\func{J}}[\mech{U}]$, and the main contribution of this paper is to bound this quantity.

\subsection{Differential Privacy}
\label{sec:differential privacy}

Differential privacy (DP) formalizes the observation that a mechanism does not reveal information about its input if the input may be replaced by a similar one without impacting the output distribution significantly. Specifically, \emph{random metric DP} \cite{Chatzikokolakis13, Dimitrakakis17, Hall12} encodes similarity via a metric on the space of input data and is a natural fit for control applications where the input comes from a metric state space:
\begin{definition}[Differential Privacy]
    Let $(\set{X}, \rho)$ be a pseudometric space and $X, \hat{X}$ be two random variables supported on $\set{X}$. A mechanism $\mech{U}: \set{X} \to \set{\Delta}(\set{U})$ is said to have \emph{$(\rho, \gamma)$-random differential privacy} ($(\rho, \gamma)$-DP) if, with probability $1 - \gamma$:
    \begin{align}
        \forall u \in \set{U}, && \log \mech{U}(X)\{u\} - \log \mech{U}(\hat{X})\{u\} \leq \rho(X, \hat{X}). \label{eq:dp}
    \end{align}
\end{definition}
Intuitively, metric DP describes a kind of Lipschitz continuity\footnote{Specifically, it is Lipschitz continuity with the specific metric on probability measures $\rho_{\mech{M}}(\mech{U}(x), \mech{U}(\hat{x})) \coloneqq \sup_{\set{V}} |\log \mech{U}(x)\{\set{V}\} - \log \mech{U}(\hat{x})\{\set{V}\}|$ and the supremum is over measurable subsets of $\set{U}$ \cite{Chatzikokolakis13}.} (with high probability) and characterizes the sensitivity of mechanisms (stochastic controllers $\mech{U}$) to replacement of $X$ with $\hat{X}$; specifically, the ratio of output densities is bounded by $\exp(\rho(X, \hat{X}))$.  An important aspect of DP is that it is preserved under composition of mechanisms~\cite{Dwork14}. Two such composition properties are:
\begin{restatable}[Post-Processing]{proposition}{postprocessing}
    If $\mech{M}_1: \set{X} \to \set{\Delta}(\set{U})$ is $(\rho, \gamma)$-DP, then for any $\mech{M}_2: \set{U} \to \set{\Delta}(\set{W})$, $\mech{M}_2 \circ \mech{M}_1$ is $(\rho, \gamma)$-DP.
\end{restatable}

\begin{restatable}[Composition]{proposition}{composition}
Let $\mech{M}_1, \mech{M}_2: \set{X} \to \set{\Delta}(\set{U})$ be $(\rho_1, \gamma_1)$-DP and $(\rho_2, \gamma_2)$-DP mechanisms respectively. Then $\mech{M} = (\mech{M}_1, \mech{M}_2)$ is $(\rho, \gamma)$-DP where $\rho(x, \hat{x}) \coloneqq \rho_1(x, \hat{x}) + \rho_2(x, \hat{x})$ and $\gamma \coloneqq \gamma_1 + \gamma_2$.
\end{restatable}

\begin{proof}
See \cref{app:privacy}.
\end{proof}

These results allow for the recursive composition of system dynamics and private controllers to yield a mechanism that computes costs in a private manner; this will allow us to extend the robustness theorem proven for single-step decision-making in the next subsection to multi-step problems (\cref{sec:multistep}).

\subsection{The Exponential Mechanism and Bounded Rationality}

A popular mechanism for DP is the \emph{exponential mechanism} \cite{Dwork14}, which provides privacy in problems where the solution of an optimization problem is desired as output. Thus, it is a natural starting point for designing a \emph{differentially-private optimal control algorithm}. Specifically, the exponential mechanism $\mech{U}^\beta: \set{X} \to \set{\Delta}(\set{U})$ is defined as the solution to the optimization problem:
\begin{align}
    \label{eq:opt_exp}
    \tag{BR-SS}
    \min_{\mech{U}}\ \func{E}\left[H(X, U)\right]\ \ \mathrm{s.t.}\ \ \func{E}[\func{D}[\mech{U}(X)||\ol{\mech{U}}]] \leq d.\nonumber
\end{align}
Here $\mech{U}^\perp \in \set{\Delta}(\set{U})$ is a ``prior" supported on $\set{U}$ independent of $X$ (emphasized by the superscript $\perp$) and $\inv{\beta}$ is the Lagrange multiplier of the relative entropy constraint.

\begin{remark}
\label{rem:bounded rationality}
This is an instance of a \emph{maximum-entropy problem} from information theory \cite{Cover99, Raginsky13, Maurer12}, and corresponds to finding a mechanism that minimizes the expected cost without significant deviation from the prior (measured by the relative entropy). As $\beta \to \infty$, $\mech{U}^\beta(x)$ simply minimizes the expected cost. As $\beta \to 0$, the solution becomes the prior. Therefore, the mechanism can be viewed as an agent employing a \emph{bounded-rational controller}, where it only uses a finite amount of information (due to computational or sensing constraints) about the state $X$ to deviate from its default behavior specified by $\ol{\mech{U}}$. Information usage (rationality) is directly controlled by $\beta$ (inversely related to $d$).
\end{remark}

The first-order optimality conditions for the problem \cref{eq:opt_exp} imply that the solution is a Gibbs measure \cite{Cover99, Raginsky13, Maurer12}:
\begin{align}
    \mech{U}^\beta(x)\{u\} &= \ol{\mech{U}}\{u\}\exp(-\beta H(x, u))\ /\ Z^\beta(x), \label{eq:fonc}\\
    Z^\beta(x) &\coloneqq \ol{\func{E}}[\exp(-\beta H(x, U))].\nonumber
\end{align}
Here, $\ol{\func{E}}[\cdot]$ is the expectation computed using $\ol{\mech{U}}$ as the controller. In the language from statistical mechanics, $\beta$ is the \emph{inverse temperature}, $H(x, u)$ is the \emph{Hamiltonian} of the system, $Z^\beta(x)$ is the \emph{partition function}, and,
\begin{align}
    F^\beta(x) \coloneqq -\inv{\beta} \log Z^\beta(x), \label{eq:free energy}
\end{align}
is known as the \emph{free energy}. The latter is notably important for a number of reasons including it being a cumulant-generating function and its equality to the Lagrangian of \cref{eq:opt_exp} conditioned on $X = x$ (see \cref{app:prop of fe}):
\begin{align}
    F^\beta(x) = \func{E}[H(x, U)] + \inv{\beta} \func{D}[\mech{U}^\beta(x)||\ol{\mech{U}}].
\end{align}

\begin{remark}
The maximum-entropy problem \cref{eq:opt_exp} is a relaxation of the \emph{rate-distortion problem} that is common in bounded rationality models \cite{Pacelli19, Pacelli20, Fox12, Tishby11}. In these problems, the relative entropy used is the mutual information between $X$ and $U$, which corresponds to an added constraint on the marginal of $\mech{U}(x)$: $\forall\set{U}' \subset \set{U}, \ol{\mech{U}}\{\set{U}'\} = \func{E}[\mech{U}^\beta(X)\{\set{U}'\}]$. Theoretically, the subsequent analysis still applies.
\end{remark}

Due to its popularity, the privacy properties of the exponential mechanism under different assumptions on $H(x, u)$ are well-studied. Only a simple assumption of Lipschitz continuity in $x$ for each $u$ is required \cite{Dimitrakakis17}. However, this assumption will not hold for the applications of interest in \cref{sec:examples} due to $\set{U}$ being non-compact. In these cases, random DP may still achievable and suitable for the purposes of robust control:
\begin{restatable}[]{proposition}{expismdp}
\label{prop:exp is mdp}
    Consider the set of all $u \in \set{U}$ for which $H(x, u)$ is at most $l$-Lipschitz in $x$, i.e. $\set{U}(l) \coloneqq \{u \in \set{U} |\ \func{L}[x \mapsto H(x, u)] < l \}$. Then, $\mech{U}^\beta(x)$ is $(2\beta l \rho, \gamma(l))$-DP, where:
    \begin{align}
        \gamma(l) \coloneqq 1 - \func{E}\left[1_{\set{U}(l)}(U) \exp\left(-2 \beta l \rho(X, \hat{\mech{X}}(X)\right)\right].\nonumber
    \end{align}
\end{restatable}

\begin{proof}
    See \cref{app:privacy}.
\end{proof}

The proposition characterizes the trade-offs in selecting the free parameter $l$ and prior $\ol{\mech{U}}$. A larger value of $l$ implies a larger set $\set{U}(l)$ but a smaller region of integration in the definition of $\gamma(l)$. For very large values of $l$, the latter will dominate and the probability that privacy fails, which is $\gamma(l)$, becomes almost certain. Note that larger $l$ implies a looser privacy constraint, since $2\beta l \rho(x, \hat{x})$ will grow for any fixed pair $x, \hat{x} \in \set{X}$. The choice of prior $\ol{\mech{U}}$ may also bias $U$ toward regions that yield a smaller Lipschitz constant.

\subsection{Quantifying the Robustness of Bounded Rationality}

The key theoretical contribution of this paper is realizing that the proposed definition of random DP can quantify the performance of the bounded-rational agent (i.e., the private controller $\mech{U}^\beta$) when only a noisy estimate $\hat{X}$ of the state $X$ is available --- thereby approximating the separation principle. This idea is combined with a large deviations argument that exploits the similarity between the large-deviation rate function \cite{Touchette16} and the free energy \cref{eq:free energy} to derive the theorem:
\begin{restatable}[]{theorem}{expisrobust}
\label{thm:single-expisrobust}
Define $\rho_{\beta}(x, \hat{x}) \coloneqq 2\beta l \rho(x, \hat{x})$. With probability at least $1 - \gamma(l)$:
\begin{align}
    \cramped{\Delta \func{J}[\mech{U}^\beta] \leq \inv{\beta} \func{E}\left[\exp\left(\rho_\beta(X, \hat{X}) + \func{D}[\mech{U}^\beta(X)||\ol{\mech{U}}]\right)\right].}\nonumber
\end{align}
\end{restatable}

\begin{proof}
See \cref{app:expisrobust}.
\end{proof}

\begin{remark}
The expectations in the theorem are implicitly conditioned on the event $U \in \set{U}(l)$, which is why a statement on expectations is shown to occur with probability $\gamma(l)$.
\end{remark}

This theorem bounds the gap $\Delta \func{J}[\mech{U}^\beta] \coloneqq \on{\func{J}}[\mech{U}^\beta] - \off{\func{J}}[\mech{U}^\beta]$ in performance between the offline and online problems in terms of $\beta$, which describes the information usage of the controller. Importantly, it relates three separate quantities: the offline expected cost (which appears in $\Delta \func{J}[\mech{U}^\beta]$), information usage (set by $\beta$), and the quality of the state estimator as measured by $\rho(x, \hat{x})$. Critically, all terms depend \emph{only on offline information} --- that is, the state estimator is not used for feedback in any of these terms but only to measure its error. Therefore, both elements of the feedback system may be designed and evaluated independently through this bound (similar to methodologies that adopt the separation principle). Moreover, the bound can be optimized to find $\beta^\star$ that yields the tightest bound on $\on{\func{J}}[\mech{U}^\beta]$.

\begin{remark}
In cases where $\func{D}[\mech{U}^\beta(X)||\ol{\mech{U}}]$ cannot be evaluated, tractable bounds may be available depending on the context (e.g. Gibbs inequality \cite{Cover99}, log-Sobolev inequalities \cite{Maurer12, Raginsky13}, and variational methods \cite{Poole19}).
\end{remark}

\section{Extension to Multi-Step Problems}
\label{sec:multistep}

In the multi-step optimal control problem, the agent attempts to minimize a sequence of non-negative cost functions $c_0, \dots, c_{t_f - 1}: \set{X} \times \set{U} \to \set{R}_+, c_{t_f}: \set{X} \to \set{R}_+$ over a time horizon $t_f$ by selecting a feedback controller $U_t \sim \mech{U}_t(X_t)$ that accounts for the system dynamics, $X_{t + 1} \sim \mech{F}_t(X_t, U_t)$. The state-input trajectories of the system generated by the choice of controller $\mech{U}_{0:t_f}$ is denoted $\mech{T}[\mech{U}_{0:t_f}]$ with $\ol{\mech{T}} \coloneqq \mech{T}[\ol{\mech{U}}_{0:t_f}]$ being the prior trajectory distribution.

Introduce the shorthand for trajectory cost:
\begin{align*}
\cramped{c_{0:t_f}(x_{0:t_f}, u_{0:t_f}) \coloneqq c_0(x_0, u_0) + \dots + c_{t_f}(x_{t_f}).}
\end{align*}
Adapting the single-step problem notation, the offline and online problems are written,
\begin{align*}%
    \min_{\mech{U}_{0:t_f}} \off{\func{J}}[\mech{U}_{0:t_f}] &\coloneqq \func{E}[c_{0:t_f}(X_{0:t_f}, U_{0:t_f})],\\
    \min_{\mech{U}_{0:t_f}} \on{\func{J}}[\mech{U}_{0:t_f}] &\coloneqq \off{\func{J}}[(\mech{U} \circ \hat{\mech{X}})_{0:t_f}],
\end{align*}
where $\hat{\mech{X}}_{0:t_f}$ are mechanisms that introduce estimation error.

Bounded rationality is achieved in a similar manner to \cref{eq:opt_exp} by constraining the relative entropy between the sequence of closed-loop and open-loop controllers:
\begin{align}
    \tag{BR-MS}
    \label{eq:brprob}
    \cramped{\min_{\mech{U}_{0:t_f}}\ \off{\func{J}}[\mech{U}_{0:t_f}]\ \ \mathrm{s.t.}\ \ \func{D}[\mech{T}[\mech{U}_{0:t_f}] || \ol{\mech{T}}] \leq d_{\mathrm{traj}}.}
\end{align}
This optimal control problem admits a recursive solution,\footnote{See, e.g. \cite{Pacelli19, Fox12, Levine18, Lambert20, Shafieepoorfard16}, for detailed solutions to similar problems.}
\begin{align}
    \mech{U}^\beta_t(x) \coloneqq \mathop{\arg \min}\limits_{\mech{U}_t}\ \func{E}\left[H_t(X, U)\right] + \inv{\beta} \func{D}[\mech{U}(x)||\ol{\mech{U}}_t],\nonumber
\end{align}
where $\beta$ is the Lagrange multiplier corresponding to the entropy constraint. Then, $\mech{U}^\beta_t$ is an exponential mechanism and both the Hamiltonian and value function $V_t(x)$ are given by the recursive equations,
\begin{align*}
    H_t(x, u) &\coloneqq c_t(x, u) + \func{E}[V_{t + 1}(\mech{F}_t(x, u))],\\
    V_t(x) &\coloneqq \func{E}[H_t(x, U_t)],
\end{align*}
where $V_{t_f}(x) \coloneqq c_{t_f}(x)$. Notably, the free energy is equivalent to the value function: $F^\beta_t(x) = V_t(x)$. The aforementioned DP composition properties allow for extension of \cref{prop:exp is mdp}, and subsequently \cref{thm:single-expisrobust}, to the multi-step problem with changes made \emph{mutatis mutandis}.
\begin{restatable}[]{proposition}{multistep}
    \label{prop:multistep}
    Let $\cramped{\set{U}_t(l_t) \coloneqq \{u\ |\ \func{L}[x\! \mapsto\! H_t (x, u)] < l_t\}}$. The mechanism,
    \begin{align}
        \mech{M}(x_{0:t_f}) \coloneqq (\mech{U}^\beta_1(x_1), \dots, \mech{U}^\beta_{t_f - 1}(x_{t_f - 1})),\nonumber
    \end{align}
    is $(\rho_\beta, \gamma)$-DP where,
    \begin{align*}
        \cramped{\rho_\beta(x_{0:t_f}, \hat{x}_{0:t_f}) \coloneqq \sum_{t = 0}^{t_f - 1} 2 \beta l_t \rho(x_t, \hat{x}_t)}, && \gamma \cramped{\coloneqq\! \sum_{t = 0}^{t_f - 1} \gamma_t,}
    \end{align*}
    and,
    \begin{align}
        \gamma_t \coloneqq 1 - \func{E}\left[1_{\set{U}_t(l_t)}(U_t) \exp\left(-2 \beta l_t \rho(X_t, \hat{\mech{X}}_t(X_t)\right)\right].\nonumber
    \end{align}
\end{restatable}

\begin{restatable}[]{theorem}{multirobust}
\label{thm:multi-expisrobust}
Let $\rho_{\beta}(x, \hat{x})$ and $\gamma$ be as in \cref{prop:multistep} and $\mech{T}^\beta \coloneqq [\mech{U}_{0:t_f}]$. With probability at least $1 - \gamma(l_t)$,
\begin{align*}
    \cramped{\Delta \func{J}[\mech{U}^\beta_{0:t_f}] \leq \frac{1}{\beta}\ \func{E}\bigg[\exp\bigg(\rho_\beta(X_{0:t_f}, \hat{X}_{0:t_f}) + \func{D}[\mech{T}^\beta || \ol{\mech{T}}]\bigg)\bigg].}
\end{align*}
\end{restatable}
\begin{proof}
See \cref{app:privacy} and \cref{app:expisrobust}.
\end{proof}
This result provides a bound on the cumulative online cost in terms of the offline cost, the level of information usage (as set by $\beta$), and the state estimation error (as measured by $\rho$) in exactly the same way as \cref{thm:single-expisrobust}.

\begin{figure*}[!t]
    \centering
    \begin{minipage}[t]{0.98\textwidth}
        \centering
        \includegraphics[width=0.48\textwidth]{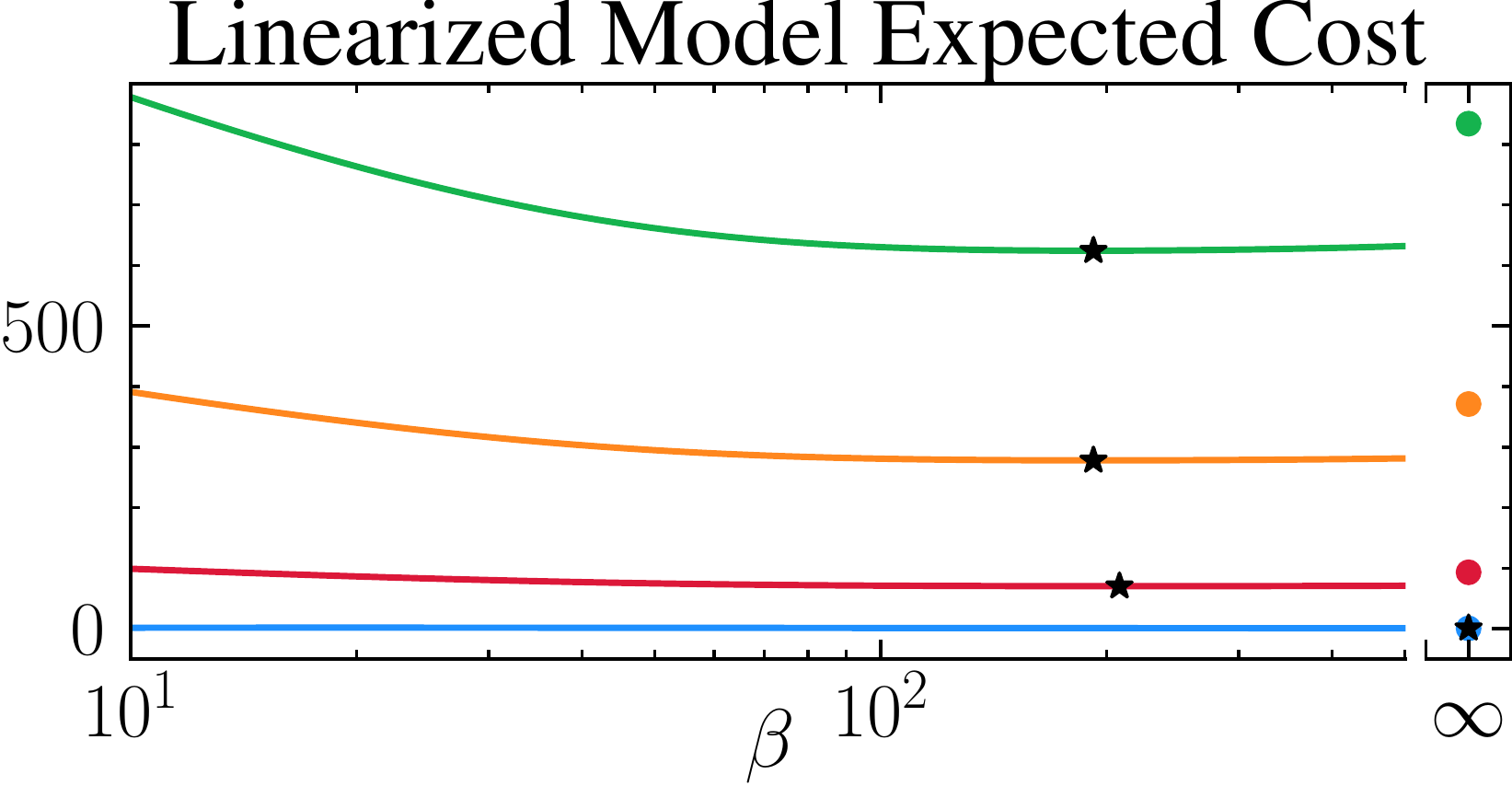}
        \hfill
        \includegraphics[width=0.48\textwidth]{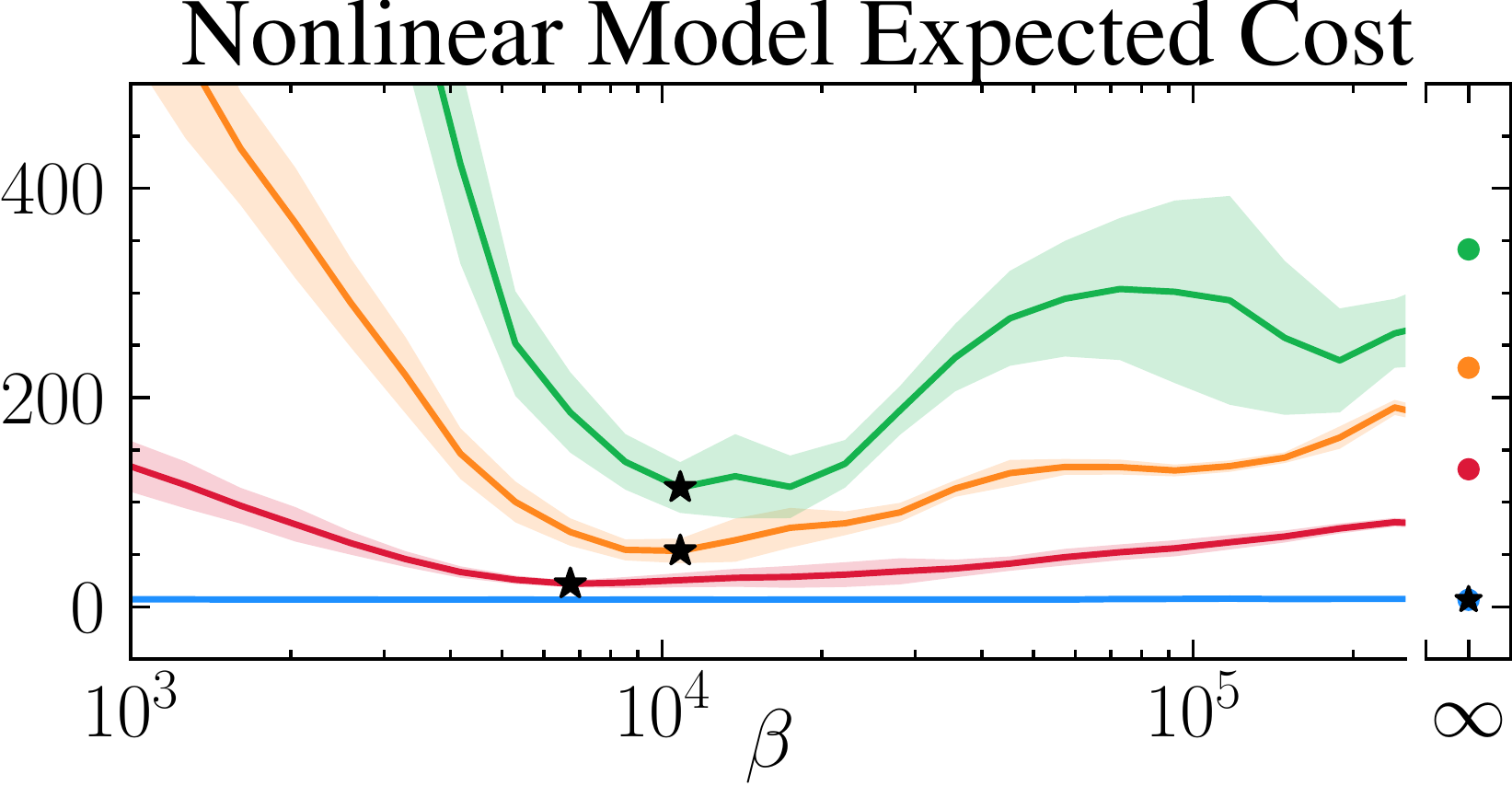}\\
        \vspace{-2mm}
        \includegraphics[width=0.5\textwidth]{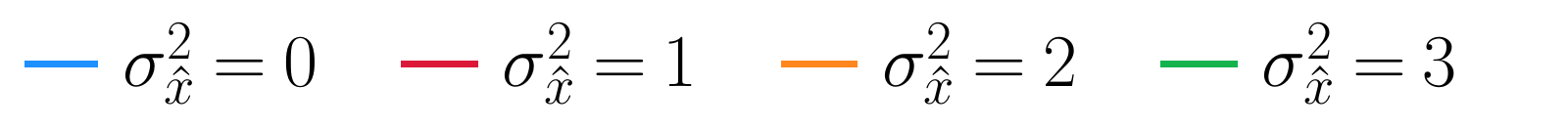}
    \end{minipage}
    \caption{Results of the numerical experiments described in \cref{sec:examples}.  The expected cost of optimally stabilizing a linearized and nonlinear quadrotor model for a range of $\beta$ and $\beta = \infty$. The latter case corresponds to the performance of the LQR and SV-MPC controllers, which are optimal given perfect state estimation. The $\star$ marker indicates $\beta^\star$ found for each value of $\sigma_{\hat{x}}^2 > 0$ within this range except for $\sigma_{\hat{x}}^2 = 0$, where the best performance is at $\beta = \infty$. For the nonlinear example, shaded regions quantify the standard deviation of the costs for 100 trials of SV-MPC. Notably, the bounded rationality controllers outperform LQR / SV-MPC in the presence of estimation error.}
    \label{fig:results}
    \vspace{-4.0mm}
\end{figure*}

        
    

\section{Numerical Examples}
\label{sec:examples}
To demonstrate the efficacy of the robustness results described in the preceding sections, three numerical examples are presented in this section: a motion planning ``double-slit'' experiment and optimal control of both the linearized and nonlinear versions of a planar quadrotor.\footnote{The source code implementing these examples using JAX \cite{Bradbury18} and all experimental parameters are listed in a publicly available repository: \href{https://github.com/irom-lab/br-dp-robust}{\texttt{https://github.com/irom-lab/br-dp-robust}}.} In each case, the dynamics of the robot are chosen to be deterministic. While the developed theory permits stochastic dynamical models, using deterministic dynamics emphasizes the sources of randomness with which this article is primarily concerned: measurement error and bounded-rational control policies.

\subsection{Motion Planning}

\textbf{Scenario.} In this problem, the objective of the robot is to find the shortest possible path to some goal region $\set{X}_g \subset \set{X}$ while avoiding a set of obstacle configurations $\set{X}_{o} \subset \set{X}$ that is disjoint from $\set{X}_g$. Both $\set{X}_g$ and $\set{X}_o$ are absorbing sets. For simplicity, the dynamics are the single integrator in the plane: $x_{t + 1} = x_t + u_t$ with $\set{X}, \set{U} = [0, 1]^2$. The problem cost is the path length if the robot reaches the goal before intersecting with an obstacle and $\infty$ otherwise. The prior $\ol{\mech{U}}_{0:t_f}$ is chosen to be a zero-mean Gaussian distribution.

The specific work space navigated by the robot is shown in \cref{fig:anchor}. The region $\set{X}_o$ consists of the boundary of the shown space and a divider punctured with a large and small slit. A path through the large slit is more robust to measurement error while a path through the smaller one is less forgiving. The goal region is in the lower right. Similar experiments are present in the literature \cite{Pedram21, Kappen05, Kappen05b, Horowitz14}.

\textbf{Results.} Importance sampling \cite{Murphy12} is used to implement the bounded rationality controller. As shown in \cref{fig:anchor}, as the information constraint tightens ($\beta \to 0$), the agent shifts from navigating the direct-but-treacherous gap to the robust (indirect) route. This demonstrates the robustness benefits of our differentially private control scheme.

\subsection{Planar Quadrotor Stabilization Problems}
\label{ssec:stabilization}

\textbf{Scenario.} These examples focus on the stabilization of a planar quadrotor (constrained to the Y-Z plane and rotating about the X axis). The system state $x \in \set{R}^6$ consists of the position and rotation of the robot and the time derivatives of these variables, i.e. $x = (y, z, \theta, \dot{y}, \dot{z}, \dot{\theta})$, while the control input $u \in \set{R}^2$ is the thrust exerted by two opposing pairs of rotors (see \cite[Section 6.6]{Steinhardt12} for dynamics). The cost functions are the linear-quadratic regulator (LQR) cost functions \cite{Anderson07}. By choosing the metric,
\begin{align*}
\rho(x, \hat{x}) = \frac{1}{2} \fronorm{x \transp{x} -  \hat{x} \transp{\hat{x}}} + \norm[2]{x - \hat{x}},
\end{align*}
the Hamiltonian $H_t(x, u)$ satisfies the conditions for the theory from \cref{sec:multistep} to apply to both systems.

The system parameters are chosen to align with the Crazyflie 2.0 quadrotors: the mass is $\SI{0.03}{kg}$ and the moment of inertia is $\SI{1.43e-5}{\kilo\gram\square\metre}$ \cite{Forster15}. 
The system is temporally discretized with a time step of $\Delta t = \SI{0.3}{s}$ and $t_f = 13$. The initial condition is sampled from small perturbations about the hover state given by the Gaussian distribution $\mech{N}(\bar{x}_0, \sigma^2_{x})$ where $\bar{y}_0 = \SI{1}{m}, \bar{z}_0 = \SI{-1}{m}$, the remaining mean entries are zero, and $\cramped{\sigma^2_{x_0} = (10^{-2}, 10^{-2}, 10^{-6}, 10^{-4}, 10^{-4}, 10^{-8})}$. The open-loop prior $\ol{\mech{U}}_{0:t_f}$ is chosen by projecting the initial distribution through the LQR solution to find the marginal input. The estimation error is drawn from a stationary Gaussian distribution. Specifically, $\hat{\mech{X}}(x) = \mech{N}(x, \sigma^2_{\hat{x}}v)$, where $\sigma^2_{\hat{x}} \in \set{R}_+$ is a scaling factor that is varied and $v = (0.25, 0.25, 0.1, 0.25, 0.25, 0.1)$.

\textbf{Linearized Results.} The system is linearized and the dynamic programming equations specifying the control policy are solved exactly for a feedback policy that is linear in the state with additive Gaussian noise. Details are included in \cref{app:lqgsln}.  The performance of the bounded rationality controller is evaluated on the linearized system and compared with a controller that is optimal assuming perfect state estimation (LQR). As shown in \cref{fig:results}, the bounded rationality controller is more robust to estimation error than the LQR controller. Moreover, the optimal information usage $\beta^\star$ decreases with the uncertainty $\sigma^2_{\hat{x}}$.

\textbf{Nonlinear Results.} The experiments were then repeated using the nonlinear planar quadrotor dynamics. There is no closed-form solution to \cref{eq:brprob}, but numerical methods for sampling from $\mech{U}^\beta_t(x)$ are available (see \cref{ssec:optimal control}). The SV-MPC algorithm \cite{Lambert20} is chosen due to its increased efficiency compared to importance sampling and the fact that it reduces to gradient-based optimization of the trajectory with multiple random initializations in the $\beta \to \infty$ limit, which is a common MPC algorithm. The results are similar to the linear case: there is a prominent local minimum for $\beta$ that outperforms MPC (which is optimal assuming perfect estimation) in the presence of estimation error when averaged over 100 trials.  In this case, the monotonic relationship between $\beta^\star$ and $\sigma^2_{\hat{x}}$ is not seen --- possibly due to SV-MPC only approximately sampling from the controller's distribution.

\section{Conclusion}
\label{sec:conclusion}
This paper proposes a new theoretical justification for the robustness of bounded-rational control policies using the framework of differential privacy. The stated performance guarantee for such policies has a modular structure reminiscent of the separation principle. Finally, multiple numerical simulations demonstrate that using differential privacy to create controllers provides robustness to estimation error.

\textbf{Future Work.} There remain a number of useful properties that need to be determined about the bounds stated in \cref{thm:single-expisrobust} and \cref{thm:multi-expisrobust} that extend beyond the scope of this paper, e.g., the tightness of the bounds and whether the relationship between $\beta^\star$ and $\sigma_{\hat{x}}^2$ is monotonic as suggested by the results for the linearized system in \cref{fig:results}. Developing a tractable method to evaluate the bounds using the methods mentioned in \cref{sec:robust} is also of great interest. Finally, there are clear opportunities to extend the experimental results of the article to new applications, such as sim-to-real transfer, and adapting other mechanisms from the DP literature to robust controller design.

\FloatBarrier

\IEEEtriggeratref{56}
\printbibliography

\newpage
\onecolumn

\setlength{\abovedisplayskip}{7pt}
\setlength{\belowdisplayskip}{7pt}

\begin{appendices}
\section{Proof of Differential Privacy Propositions}
\label{app:privacy}
\postprocessing*

\begin{proof}
    Let $\set{W}' \subset \set{W}$ be measurable and $\set{U}(\set{W}') \coloneqq \{u \in \set{U}\ |\ \mech{M}_2(u)\{\set{W}'\} > 0\}$. Assume momentarily that $\mech{M}_1$ is always private, i.e. $\gamma = 1$. Then, for any $x, \hat{x} \in \set{X}$ and $U \coloneqq \mech{M}_1(x), \hat{U} \coloneqq \mech{M}_2(\hat{x})$,
    \begin{align*}
        \frac{(\mech{M}_2 \circ \mech{M}_1)(x)\{\set{W}'\}}{(\mech{M}_2 \circ \mech{M}_1)(\hat{x})\{\set{W}'\}} &= \frac{\mathbb{E}[1_{\set{U}(\set{W}')}(U) \cdot \mech{M}_1(x)\{U\} \cdot \mech{M}_2(U)\{\set{W}'\}]}{\mathbb{E}[1_{\set{U}(\set{W}')}(U) \cdot \mech{M}_1(\hat{x})\{U\} \cdot \mech{M}_2(U)\{\set{W}'\}]},\\
        &\leq \frac{\mathbb{E}[1_{\set{U}(\set{W}')}(U) \cdot \exp(\rho(x, \hat{x})) \mech{M}_1(\hat{x})\{U\} \cdot \mech{M}_2(U)\{\set{W}'\}]}{\mathbb{E}[1_{\set{U}(\set{W}')}(U) \cdot \mech{M}_1(\hat{x})\{U\} \cdot \mech{M}_2(U)\{\set{W}'\}]},\\
        &\leq \exp(\rho(x, \hat{x})) \frac{\mathbb{E}[1_{\set{U}(\set{W}')}(U) \cdot \mech{M}_1(\hat{x})\{U\} \cdot \mech{M}_2(U)\{\set{W}'\}]}{\mathbb{E}[1_{\set{U}(\set{W}')}(U) \cdot \mech{M}_1(\hat{x})\{U\} \cdot \mech{M}_2(U)\{\set{W}'\}]} = \exp(\rho(x, \hat{x})).
    \end{align*}
    When $\gamma < 1$, the same argument holds probability $1 - \gamma$ and $x, \hat{x}$ replaced by random variables.
\end{proof}

\composition*

\begin{proof}
    Assume $\gamma_1, \gamma_2 = 1$. Let $\set{V} \coloneqq \set{U}_1 \times \set{U}_2 \subset \set{U} \times \set{U}$ be measurable and $x, \hat{x} \in \set{X}$. Then:
    \begin{align*}
        \log \frac{\mech{M}(x)\{\set{U}\}}{\mech{M}(\hat{x})\{\set{U}\}} &= \log \frac{\mech{M}_1(x)\{\set{U}_1\}}{\mech{M}_1(\hat{x})\{\set{U}_1\}} + \log \frac{\mech{M}_2(x)\{\set{U}_2\}}{\mech{M}_2(\hat{x})\{\set{U}_2\}} \leq \rho_1(x, \hat{x}) + \rho_2(x, \hat{x}).
    \end{align*}
    Now, when $\gamma_1, \gamma_2 < 1$, the probability the above holds is estimated via the union bound:
    \begin{align*}
        \func{P}\{\mech{M}_1\ \textrm{is private}\ \wedge\ \mech{M}_2\ \textrm{is private}\} &= 1 - \func{P}\{\mech{M}_1\ \textrm{is not private}\ \vee\ \mech{M}_2\ \textrm{is not private}\},\\
        &\geq 1 - (\func{P}\{\mech{M}_1\ \textrm{is not private}\} + \func{P}\{\mech{M}_2\ \textrm{is not private}\}) \geq 1 - (\gamma_1 + \gamma_2).
    \end{align*}
\end{proof}

\expismdp*

\begin{proof}
For all $u \in \set{U}(l)$, \cref{eq:dp} holds for the metric $2 \beta l \rho(x, x')$ \cite[Theorem 8]{Dimitrakakis17}. All that is needed is to lower bound the probability of sampling $U \in \set{U}(l)$, i.e. $U \in \mech{U}(\hat{X})\{\set{U}(l)\}$. For any $x, \hat{x} \in \set{X}$:
\begin{align*}
    \mech{U}^\beta(x)\{\set{U}(l)\} &\leq \exp(2 \beta l \rho(x, \hat{x})) \mech{U}^\beta(\hat{x})\{\set{U}(l)\}
\end{align*}
This inequality may be rewritten as:
\begin{align*}
    \ol{\func{E}}\left[1_{\set{U}(l)}(U) \exp\left(-2 \beta l \rho(x, \hat{x})\right) \frac{\exp(-\beta H(x, U)}{Z^\beta(x)}\right] &\leq \ol{\func{E}}\left[1_{\set{U}(l)}(U) \frac{\exp(-\beta H(\hat{x}, U)}{Z^\beta(\hat{x})}\right].
\end{align*}
Apply the law of total expectation:
\begin{align*}
    1 - \gamma(l) = \ol{\func{E}}\left[1_{\set{U}(l)}(U) \exp\left(-2 \beta l \rho(X, \hat{X})\right) \frac{\exp(-\beta H(X, U)}{Z^\beta(X)}\right] &\leq \ol{\func{E}}\left[1_{\set{U}(l)}(U) \frac{\exp(-\beta H(\hat{X}, U)}{Z^\beta(\hat{X})}\right].
\end{align*}
The right hand side of this inequality is the probability that $U \in \set{U}(l)$ as desired.
\end{proof}

\multistep*

\begin{proof}
Recursively apply the composition property and the fact that $\mech{U}^\beta$ is $(2\beta l_t \rho, \gamma(l_t))$-DP.
\end{proof}


\section{Variational Representation of the Free Energy}
\label{app:prop of fe}
\begin{proposition}
\label{prop:free energy variational}
For any $\mech{U}(x)$ such that for all $x \in \set{X}$, $\mech{U}(x)$ is absolutely continuous with respect to $\ol{\mech{U}}$ (denoted $\mech{U}(x) \ll \ol{\mech{U}}$), it holds that:
\begin{align}
\forall x \in \statespace, u \in \inputspace, && F^\beta(x) \leq \func{E}[H(x, \mech{U}(x))] + \frac{1}{\beta} \func{D}[\mech{U}(x)||\ol{\mech{U}}]. \label{eq:free energy variational inequality}
\end{align}
This inequality is tight when $\mech{U} = \mech{U}^\beta$.
\end{proposition}

\begin{proof}
     To begin, rearrange \cref{eq:fonc} into:
     \begin{align*}
         H(x, u) + \frac{1}{\beta} \log \frac{\mech{U}^\beta(x)\{u\}}{\ol{\mech{U}}\{u\}} &= F^\beta(x).
     \end{align*}
     Take the expectation with respect to $\mech{U}^\beta(x)$ of both sides ($F^\beta(x)$ is constant in $u$) to achieve the equality half of the consequent \cref{eq:free energy variational inequality}. Now, since $\mech{U}^\beta$ is the solution to the minimization problem \cref{eq:opt_exp}, any other feasible choice of $\mech{U}$ (i.e. any choice such that $\mech{U} \ll \ol{\mech{U}}$) will produce a value for the problem greater than or equal the optimal value. Therefore,
     \begin{align*}
         F^\beta(x) &= \func{E}[H(x, \mech{U}^\beta(x))] + \inv{\beta} \func{D}[\mech{U}^\beta(x)||\ol{\mech{U}}]\\
         &\leq \func{E}[H(x, \mech{U}(x))] + \frac{1}{\beta} \DD[\mech{U}(x)||\ol{\mech{U}}].
     \end{align*}
\end{proof}

\section{Proof of Robustness Theorems}
\label{app:expisrobust}
\expisrobust*
\begin{proof}
Define the standard notation $\func{P}\{\set{A}\} \coloneqq \func{E}[1_{\set{A}}(\set{A})]$ for $\set{A}$ is a subset of the space of states, state estimates, and inputs, i.e., $\set{A} \subseteq \set{X} \times \set{X} \times \set{U}$. The mechanisms,
\begin{align*}
    \mech{H}(\hat{x}; x) \coloneqq H\left(x, \mech{U}^\beta(\hat{x})\right), && \mech{H}(x) \coloneqq \mech{H}(x; x), && \hat{\mech{H}}(x) \coloneqq \mech{H}\big(\mech{\hat{X}}(x); x\big).
\end{align*}
are the cost incurred by a particular $(\hat{x}, x)$ pair, the offline cost, and the online cost respectively. For any $x \in \set{X}$, the mechanism $\hat{x} \mapsto \mech{H}(\hat{x}; x)$ is $(2\beta l \rho, \gamma(l))$-DP due to the post-processing property (\cref{app:privacy}). As a result, for any state $x \in \set{X}$ and any measurable set of costs $\set{C} \subset \set{R}_+$,
\begin{align}
     \log\frac{\hat{\mech{H}}(x)\{\set{C}\}}{\mech{H}(x)\{\set{C}\}} \leq \rho_\beta(x, \hat{X}), \label{eq:change of measure}
\end{align}
holds with probability $1 - \gamma(l)$. That is, with high probability, the measure over costs online costs is exchanged with the measure over offline costs with a bounded change in outcome probabilities.

Now, the Chernoff bound \cite{Raginsky13} applied to the mechanism $\hat{\mech{H}}(x)$ states:
\begin{align*}
     \forall a, \zeta \in \set{R}_+, && \func{P}\left\{\hat{\mech{H}}(x) \geq a \right\} \leq \func{E}\left[\exp\left(\zeta \hat{\mech{H}}(x) - \zeta a\right)\right].
\end{align*}
Apply \cref{eq:change of measure} to switch to the offline setting with probability $1 - \gamma(l)$:
\begin{align*}
     \func{P}\left\{\hat{\mech{H}}(x) \geq a \right\} \leq \func{E}\left[\exp\left(\rho_\beta(x, \hat{X}) + \zeta \mech{H}(x) - \zeta a\right)\right].
\end{align*}
Since $\rho_\beta(x, \hat{X})$ is independent from $\mech{H}(x)$, the expectation factors and the right-hand side may be rewritten,
\begin{align}
    \func{P}\left\{\hat{\mech{H}}(x) \geq a \right\} &\leq \func{E}\left[\exp\left(\rho_\beta(x, \hat{X}) - \zeta a\right)\right]\func{E}\left[\exp\big(\zeta \mech{H}(x)\big)\right],\nonumber\\
    &= \func{E}\left[\exp\left(\rho_\beta(x, \hat{X}) + I(a; \zeta)\right)\right],  \label{eq:concentration}
\end{align}
where:
\begin{align*}
    I(a; \zeta) \coloneqq \log \func{E}\left[\exp\bigg(\zeta \mech{H}(x)\bigg)\right] - \zeta a.
\end{align*}

When the infimum of $I(a; \zeta)$ is taken with respect to the $\zeta$, this function is known as the \emph{rate function} in large-deviation theory \cite{Touchette16}. However, such extremization is not performed because the minimizing $\gamma$ to find the minimizing $\gamma^\star$ since it will be dependent on both $x$ and $\beta$ and is difficult to analyze. Instead, the values of $a$ and $\zeta$ are chosen based on the properties of $\mech{U}^\beta(x)$ so that $\beta$ may be optimized over to jointly minimize all terms. In addition to be analytically tractable, it prevents $\beta^\star$ being tied to $x$ through $\zeta^\star$. If $\beta^\star$ is dependent on $x$, the mechanism $\mech{U}^{\beta^\star(x)}(x)$ may not in fact be differentially private. Specifically, choose $\zeta = \beta$ and $a = \func{E}[\mech{H}(x)] + b$. In this case, as a result of \cref{prop:free energy variational},
\begin{align}
    I\left(\func{E}[\mech{H}(x)] + b; \beta \right) &\coloneqq \log \func{E}\left[\exp\big(\beta H(x, U)\big)\right] - \beta (\func{E}[\mech{H}(x)] + b),\nonumber\\
    &= \log \ol{\func{E}}\left[\frac{\exp\big(\beta H(x, U)\big) \exp\big(-\beta H(x, U)\big)}{Z^\beta(x)}\right] - \beta (\func{E}[\mech{H}(x)] + b),\nonumber\\
    &= -\log Z^\beta(x) - \beta (\func{E}[\mech{H}(x)] + b),\nonumber\\
    &= \beta F^\beta(x) - \beta (\func{E}[\mech{H}(x)] + b),\label{eq:change 1}\\
    &= \func{D}[\mech{U}^\beta(x) || \ol{\mech{U}}] - \beta b.\label{eq:change 2}
\end{align}
Plugging these values into  \cref{eq:concentration} yields:
\begin{align*}
    \func{P}\left\{\hat{\mech{H}}(x) - \func{E}[\mech{H}(x)] \geq b \right\} &\leq \func{E}\left[\exp\left(\rho_\beta(x, \hat{X}) + \func{D}[\mech{U}^\beta(x) || \ol{\mech{U}}] - \beta b\right)\right].
\end{align*}
After integrating over $b$ from $0$ to $\infty$ to find,
\begin{align*}
    \func{E}[\hat{\mech{H}}(x)] - \func{E}[\mech{H}(x)] &\leq \inv{\beta}\func{E}\left[\exp\left(\rho_\beta(x, \hat{X}) + \func{D}[\mech{U}^\beta(x) || \ol{\mech{U}}]\right)\right].
\end{align*}
Apply the law of total expectation to get the desired result:
\begin{align*}
    \on{\func{J}}[\mech{U}^\beta] - \opt{\func{J}} \leq \inv{\beta} \func{E}\left[\exp\left(\rho_\beta(X, \hat{X}) + \func{D}[\mech{U}^\beta(X) || \ol{\mech{U}}]\right)\right].
\end{align*}
\end{proof}

\multirobust*

\begin{proof}
Let $U_t \sim \mech{U}^\beta_t(X_t)$, $\hat{U}_t \sim (\mech{U}_t^\beta \circ \hat{\mech{X}}_t)(X_t)$, and let $\func{P}\{\cdot\}$ be the same as in the proof of \cref{thm:single-expisrobust} except now defined for trajectories of the pertinent variables, i.e. $X_{0:t_f}, \hat{X}_{0:t_f}, U_{0:t_f}$. Similarly, define:
\begin{align*}
    \mech{H}(\hat{x}_{0:t_f}; x_{0:t_f}) \coloneqq \sum_{t = 0}^{t_f - 1} c_t\left(x_t, \mech{U}^\beta_t(\hat{x}_t)\right) + c_{t_f}(x_{t_f}), && \mech{H}(x_{0:t_f}) \coloneqq \mech{H}(x_{0:t_f}; x_{0:t_f}), && \hat{\mech{H}}(x_{0:t_f}) \coloneqq \mech{H}\big(\mech{\hat{X}}_{0:t_f}(x_{0:t_f}); x_{0:t_f}\big).
\end{align*}

The exact same procedure as in \cref{thm:single-expisrobust} is now applied except with trajectories of the pertinent variables. The only step that needs to be justified is moving from \cref{eq:change 1} to \cref{eq:change 2}. However, this step is handled by using the inequality form of \cref{prop:free energy variational}. For clarity, \cref{eq:brprob} is rewritten in the following explicit form:
\begin{mini}[2]{\mech{T},\ \mech{U}_{0:t_f}}{\func{E}\left[c_{0:t_f}(X_{0:t_f}, U_{0:t_f})\right]}{\tag{BR-MS}\label{eq:brprob-rewritten}}{}
    \addConstraint{X_0 \sim \mech{X}_0,}{}{X_{t + 1} \sim \mech{F}_t(X_t, U_t)}
    \addConstraint{U_t \sim \mech{U}_t(X_t),}{}{\mech{T}_t \ll \ol{\mech{T}}_t}
    \addConstraint{\mech{T}(x_0) = (x_0, X_{1:t_f}, U_{0:t_f}),\quad}{}{\func{D}[\mech{T}(x_0)||\ol{\mech{T}}]}{\leq d_{\mathrm{traj}}.}
\end{mini}
Here, the trajectory distribution, conditioned on the initial condition $x_0$, is the mechanism $\mech{T}(x_0)$ and appears as a decision variable. Through the equality constraint, it is entirely defined by the choice of $\mech{U}_{0:t_f}(x_{0:t_f})$. However, \cref{eq:brprob-rewritten} is now in the form of a maximum-entropy problem, like \cref{eq:opt_exp}, except it has an additional dynamical constraint: $X_{t + 1} \sim \mech{F}_t(X_t, U_t)$. Since any solution $\mech{T}^\star(x_0)$ to \cref{eq:brprob-rewritten} will yield an objective greater than or equal to the value of the problem without the dynamical constraint, it follows from \cref{prop:free energy variational} that:
\begin{align*}
    \log \func{E}\left[\exp(\beta c_{0:t_f}(X_{0:t_f}, U_{0:t_f}))\right] &= -\log \ol{\func{E}}\left[\exp\left(-\beta c_{0:t_f}(X_{0:t_f}, U_{0:t_f})\right)\right] \leq \beta \func{E}\left[c_{0:t_f}(X_{0:t_f}, U_{0:t_f})\right] + \func{D}[\mech{U}^\beta(x) || \ol{\mech{U}}].
\end{align*}
Therefore, the sign of \cref{eq:change 2} is now an inequality, but its direction coincides with that of \cref{eq:concentration} so the proof is not impacted.
\end{proof}

\section{Explicit LQG Solution of the Bounded-Rationality Policy}
\label{app:lqgsln}
This section derives the linear-quadratic-Gaussian (LQG) solution for the bounded-rationality optimal control problem \cref{eq:brprob}. The strict notation the document's body is relaxed for one more suited for linear algebra. For completeness the model is restated with the addition of Gaussian process noise. The dynamics are,
\begin{align}
    X_{t + 1} = A X_t + B U_t + \epsilon_t, && \epsilon_t \sim \mech{N}(0, \Sigma_{\epsilon_t}),
\end{align}
and the cost functions are,
\begin{align}
    c_t(x, u) = \transp{x} Q x + \transp{u} R u, && c_{t_f}(x) = \transp{x} Q_f x.
\end{align}
Here, $\set{X} = \set{R}^n, \set{U} = \set{R}^m$ and $Q, Q_f, \Sigma_{\epsilon_t} \in \set{S}_{++}^{n}, R \in \set{S}_{++}^{m}$ are positive-definite matrices, and $X_t, U_t$ are the random states at time $0 \leq t \leq t_f$.

The goal is to find a policy that minimizes the total of the above cost functions over the horizon $t_f$. To begin the dynamic programming argument, assume a quadratic form for $V_{t + 1}(x)$:
\begin{align}
    V_{t+1}(x) = \frac{1}{2}\transp{x} P_{t + 1} x + \transp{b}_{t + 1} x + d_{t + 1}.
\end{align}
Then, $H_t(x, u)$ is expressed as,
\begin{align}
    H_t(x, u) &= \frac{1}{2}\transp{(u - v_t)} \inv{S}_t (u - v_t) + k_t(x),
\end{align}
where:
\begin{align}
    S_t \coloneqq \inv{(\transp{B}P_{t + 1}B + R)}, && v_t \coloneqq -S_t\transp{B}(P_{t + 1}Ax + b_{t + 1}),
\end{align}
Since $H_t(x, u)$ is quadratic, choosing $\ol{\mech{U}}_t$ to be Gaussian ensures \cref{eq:fonc} is Gaussian: $\pi_t(u_t | x_t) = \mech{N}(\bar{u}_t(x_t), \Sigma_{u_t|x_t})$. Plugging in the prior and completing the square yields:
\begin{align}
    \bar{u}(x) &= -\Sigma_{u_t|x_t}(\beta \transp{B}P_{t + 1}Ax + \beta \transp{B} b_{t + 1} + \inv{\Sigma}_{\ol{u}_t} \ol{\bar{u}}_t), && \inv{\Sigma}_{u_t|x_t} = \beta \transp{B}P_{t + 1}B + \beta R - \inv{\Sigma}_{\ol{u}_t}.
\end{align}
More practically, the distribution may be expressed as the affine relationship,
\begin{align}
    u_t = K_t x_t + \eta_t, && \eta_t \sim \mathcal{N}(\bar{\eta}_t, \Sigma_{\eta_t}).
\end{align}
where:
\begin{align}
    \inv{\Sigma}_{\eta_t} &= \beta \transp{B}P_{t + 1}B + \beta R + \inv{\Sigma}_{\ol{u}_t}, && \bar{\eta}_t = -\Sigma_{\eta_t}(\beta \transp{B}b_{t + 1} - \inv{\Sigma}_{\ol{u}_t}\ol{\bar{u}}_t), && K_t =  -\beta \Sigma_{\eta_t}\transp{B}P_{t + 1}A.
\end{align}
Finally, the value function at time $t$ can be computed as,
\begin{align}
    V_t(x) = \func{E}[H_t(x, U_t)] = \frac{1}{2}\transp{x} P_t x + \transp{b}_t x + d_t.
\end{align}
where
\begin{align}
    P_t &= Q + \transp{K}_t R K_t + \transp{A} P_{t + 1} A + \transp{K}_t \transp{B} P_{t + 1} B K_t + \frac{1}{2}(\transp{K}_t \transp{B} P_{t + 1} A + \transp{A}P_{t + 1}BK_t),\nonumber\\
    \transp{b}_t &= \transp{\bar{\eta}}_t R K_t + \transp{\bar{\eta}}_t \transp{B} P_{t + 1} A + \transp{\bar{\eta}}_t \transp{B} P_{t + 1} B K_t + \transp{b}_{t + 1}(A + BK_t),\nonumber\\
    d_t &= d_{t + 1} + \frac{1}{2}\bar{\eta}_t (R + \transp{B} P_{t + 1} B) \bar{\eta}_t + \transp{b}_{t + 1} B \bar{\eta}_t + \frac{1}{2} \mathrm{tr}\left(\Sigma_{\eta_t}(R + \transp{B}P_{t + 1} B)\right) + \frac{1}{2}\mathrm{tr}(\Sigma_{\epsilon_t} P_{t + 1}) .\nonumber
\end{align}
The recursive description is completed by observing that the terminal cost function implies $P_{t_f} = Q_f, b_{t_f} = 0, d_{t_f} = 0$.
\begin{remark}
As $\beta \to \infty$, the LQR solution is recognized (and is the same as the linear-quadratic Gaussian problem solution). Similarly, as $\beta \to 0$, it is clear that the statistics of $U_t$ match those of $\ol{\pi}_t$.
\end{remark}
\end{appendices}

\end{document}